# Chance and the Continuum Hypothesis

DANIEL HOEK, PRINCETON UNIVERSITY, SEPTEMBER 2020, DANIELHOEK.COM. PLEASE CITE PUBLISHED VERSION.[1]

ABSTRACT: This paper presents and defends an argument that the continuum hypothesis is false, based on considerations about objective chance and an old theorem due to Banach and Kuratowski. More specifically, I argue that the probabilistic inductive methods standardly used in science presuppose that every proposition about the outcome of a chancy process has a certain chance between 0 and 1. I also argue in favour of the standard view that chances are countably additive. Since it is possible to randomly pick out a point on a continuum, for instance using a roulette wheel or by flipping a countable infinity of fair coins, it follows, given the axioms of ZFC, that there are many different cardinalities between countable infinity and the cardinality of the continuum.

Quine told us that there are no islands in science, and that even our mathematical beliefs have to face the tribunal of experience. The examples usually adduced in support of his view are non-Euclidean geometry and quantum probability, mathematical theories whose development was in part inspired by physical discoveries. But in those cases, it is unclear whether our beliefs about the mathematics itself were revised, or just our beliefs about the applicability of certain mathematical theories. In this paper I give a clearer demonstration of the way the empirical can impinge on the mathematical. I argue that our scientific knowledge of chancy processes provides compelling abductive reason to think that uncountable collections of points on a continuum come in a wide variety of sizes, contradicting Cantor's famous hypothesis that all of them have the same size. I will not establish any novel mathematical results: my argument is based on a theorem due originally to Stefan Banach and Casimir Kuratowski (1929), which is well known to set theorists but less well known to philosophers.

There is room for an argument of this kind because, mathematically speaking, the question of the continuum hypothesis is in a very strong sense unresolved. Gödel (1938) and Cohen (1963) proved that the hypothesis is independent of the ZFC axioms (Zermelo-Fraenkel set theory with Choice). This showed that any argument to settle the issue would have to appeal to assumptions that go beyond those canonised in standard set theory. But it remained opaque what outside source those principles might spring from. For a while some hoped, with Gödel, that internal, mathematical considerations would settle the matter in the form of a large cardinal axiom. But nowadays there is reason to be pessimistic on that front: all large cardinal axioms that have been tried leave the undecidable status of the continuum hypothesis unaffected, and there are generally applicable strategies for extending that conclusion to any further large cardinal axioms (Levy and Solovay 1967, Honzik 2017). More recently,

[1] *Philosophy and Phenomenological Research*, online first, 2020: 1-22. I would like to thank Jonathan Hickman for introducing me to Banach and Kuratowski's theorem; this paper has also benefited a great deal from conversations with and comments from Andrew Bacon, David Chalmers, Nicholas DiBella, Cian Dorr, Adam Elga, Colin Elliot, Jeremy Goodman, Zachary Goodsell, Thomas Hofweber, Adam Lovett, Lorenzo Rossi, Chris Scambler, Trevor Teitel, Dan Waxman, Snow Zhang, participants in NYU's Washington Square Circle and Rutgers' Foundations of Probability Group, and an anonymous reviewer.

Hugh Woodin has explored mathematically motivated argument strategies for refuting the continuum hypothesis (Woodin 2001) and for establishing it (Woodin 2010). But even by his own lights, these strategies have hitherto delivered inconclusive results (Rittberg 2015, Woodin 2019). In 1986, Chris Freiling pointed to a very different avenue, by suggesting that we can settle the question by appeal to considerations about the incidence of certain random events. Here I provide an argument in that tradition, which lacks the problematic features of Freiling's specific proposal.

Certain random processes have a continuum of possible outcomes. For instance, if you spin a roulette wheel, the "0" slot could end up in a continuum of different positions once the wheel has come to a standstill. A random dart may land in one of a continuum of different positions on the dartboard. And if you have an endless row of coins and flip all of them, the outcome is one of a continuum of possible sequences of heads and tails. Quantum mechanical examples of such processes include a momentum measurement on a particle with known position, or a countable infinity of $x$-spin measurements of particles in $z$-spin eigenstates. I will argue that the chance distributions of such procedures have a set of formal features that, jointly, tell us something about the size of their outcome space — and thus about the size of the continuum.

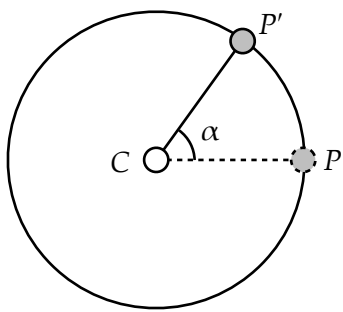

*Diagram 1: Any whirl of a spinner determines a plane angle α = ∠PCP' between the pointer's initial position P, the centre of rotation C, and the pointer's final position P'*

For perspicuity, I will focus primarily on the example of a *spinner*: a rotating device with an off-centre pointer, which on a given spin lands some random angle away from its starting point. A roulette wheel is an example, or a spinning top, or one of those fidget spinners you can buy at a gas station. Suppose we have a concrete object like this in front of us, mark a pointer, and give it a whirl. The *outcomes* of the procedure are the angles at which the pointer could land. The set Ω = [0°, 360°) of all those outcomes is the *outcome space*. *Events* are represented by subsets of Ω: for instance, the interval (10°, 15°) represents the event that the pointer lands at an angle between 10° and 15°. Standardly, events are taken to be propositions about where the spinner lands on a particular occasion (e.g. Lewis 1980). One can also think of them more abstractly as "ways the spinner might land". For our purposes that difference does not matter. All that matters is that events are the bearers of chance.

The *chance* of an event is a real number between 0 and 1 that indicates the objective likelihood that the event takes place. For instance, if our spinner has an uneven mass distribution, the chance of the event $[0°, 180°)$ might be 0.53, meaning that the spinner is objectively more likely than not to land on that side. My argument will rely on a realist view of chances, according to which they supervene on objective, physical features of the setup; for reasons discussed in Sections I and II, I take this realism about chance to be justified by the need to explain the overwhelming success of probabilistic methods of prediction in science. In particular, this view entails that chances are not "assigned" or "defined" by us any more than the mass and temperature distribution of the spinner are, and that they are not in any straightforward sense a measure of our actual expectations. The *chance distribution* is simply the function from events to real numbers that relates events to their chances.

Below I shall argue that we have good reason to think the chance distribution Ch of a random spinner has the following four properties:

A. *Perfect Precision*. It is certain that exactly one of the outcomes occurs: $\mathrm{Ch}(\Omega) = 1$.
B. *Improbable Outcomes*. No individual outcome has a positive chance of being realised: $\mathrm{Ch}(\{x\}) = 0$ for any outcome $x \in \Omega$.
C. *Totality*. Every event $E \subseteq \Omega$ has a chance $\mathrm{Ch}(E)$.
D. *Countable Additivity*. If $E_0, E_1, \ldots \subseteq \Omega$ are disjoint, then $\mathrm{Ch}(E_0 \cup E_1 \cup \ldots) = \mathrm{Ch}(E_0) + \mathrm{Ch}(E_1) + \ldots$

By far the most contentious of these claims is the totality premise (C), which I will defend in detail in Section II. Section I will cover (A) and (B), and Section III argues for premise (D).

Mathematically, the joint upshot of (A-D) is that a continuum-sized set $\Omega$ admits of a total, countably additive measure $\mathrm{Ch}: \mathscr{P}(\Omega) \to [0, 1]$ with the property that $\mathrm{Ch}(\{x\}) = 0$ for every $x \in \Omega$. Now in 1930, Ulam showed, within ZFC, that this is not true of the first uncountable cardinal $\aleph_1$: that is, $\aleph_1$-sized sets do *not* admit a measure of this kind (Jech 2002, lemma 10.13, p. 132). So if (A-D) are true, it follows that the continuum $2^{\aleph_0}$ is distinct from the first uncountable cardinal $\aleph_1$, which is just to say that the continuum hypothesis $2^{\aleph_0} = \aleph_1$ is false. In fact, it turns out to follow that $2^{\aleph_0}$ is *far* greater than $\aleph_1$, exceeding not only the second uncountable cardinal $\aleph_2$, but also $\aleph_3, \aleph_4 \ldots \aleph_\omega, \aleph_{\omega_1}, \ldots \aleph_{\omega_\omega}, \aleph_{\omega_{\omega_\omega}}, \ldots$ and the least weakly inaccessible cardinal.

Let me start with a few preliminaries about the argument. Firstly, I do *not* assume the spinner is symmetric or fair: the premises are consistent with a heavily skewed chance distribution. Secondly, as is standard, I take it that events can be physically possible and still have chance 0. After all, a whirl of a spinner has uncountably infinitely many possible outcomes; if they all had a positive chance of being realised, the sum would exceed 1. Some theorists find this counterintuitive, and have developed treatments of probability that allow for *infinitesimal* chances: non-zero chances smaller than any positive number. The chance of an event is then not a real number but rather a non-standard real or hyperreal number (Benci et al. 2013, Hofweber and Schindler 2016; for objections see Williamson 2007, Easwaran 2014). This makes no difference for the present argument: even if chances are hyperreals, they still have a *real part* — the unique real number closest to their hyperreal value. And thus, as I will argue, a real-valued measure Ch satisfying (A-D) would still have to exist — it is just that Ch($E$) would represent the *real part* of the chance of $E$ rather than the chance itself.

Thirdly, let me say a bit more about the totality premise (C), which basically says that any proposition about the outcome of the spin has a particular chance, possibly 0, of coming true. While this is a very intuitive assumption, it is liable to be controversial in the case of infinite-outcome procedures. The default mathematical treatment of chance in such cases straightforwardly incorporates premises (A), (B) and (D). But it does not bear out premise (C) in the same way: when writing down a probability distribution on a continuum-sized outcome space, the custom is to single out a class of "measurable" events first, and to assign chances only to those events. In Section II, I argue that this is done for expediency, and that there is no good reason to think there is a real distinction in the world between chance-bearing and chance-free events corresponding to the mathematical distinction between measurable and non-measurable sets of points. Furthermore, I argue that there is strong abductive reason to think (C) is true, because the assumption of totality is implicit in our ordinary scientific reasoning about chance.

Fourthly, while I presented this as an argument about a random spinner, its premises (A-D) are stated schematically, and $\Omega$ can stand in for other continuum-sized outcome spaces. For instance, (A-D) remain plausible if one interprets $\Omega$ as the outcome space of a random dart, or an infinite sequence of coin flips, or one of the quantum-mechanical processes mentioned. These latter substitutions are particularly relevant if you think that the chances of quantum-mechanics are the only genuine chances in nature: to establish the falsehood of the continuum hypothesis, only one version of the argument needs to hold up. You can also substitute a merely possible process. Suppose, say, that a spinner satisfying conditions (A-D) is metaphysically possible. Then it is at least possible for a chance measure of the appropriate kind to exist. But the existence of such a purely mathematical function is not a contingent matter. So the measure actually exists, whence the continuum hypothesis is actually false. In order to evade that version of the argument, it is not enough to deny that actual, real-world spinners fail to meet one of the conditions (A-D). One has to deny that such spinners are even possible. Likewise, infinite arrays of flipping coins, random darts, and their quantum-mechanical analogues would somehow have to be impossible.

Formally, the argument motivates the addition of the following axiom to ZFC:

> *Chance Measurability of the Continuum*: Continuum-sized sets $\Omega$ admit of a total, countably additive measure Ch such that $Ch(\{x\}) = 0$ for every $x \in \Omega$. (M)

Solovay (1971) showed that the theory ZFC + M is consistent.[2] Banach and Kuratowski (1929) first proved that the negation of the continuum hypothesis, $2^{\aleph_0} > \aleph_1$, is a theorem of ZFC + M. In an appendix I provide an intuitive rendition of their proof. Ulam (1930) and Solovay (1971) strengthened this result, showing that given ZFC + M, there must be a cardinal $\kappa \leq 2^{\aleph_0}$ such that $\kappa$ is (i) weakly inaccessible, (ii) real-valued measurable[3] and (iii) weakly Mahlo. These are *large cardinal* properties in the sense that unlike ZFC + M, ZFC alone cannot prove that any of them are instantiated. But $\kappa$ is not a large cardinal in the sense of implying the consistency of ZFC.

Since the proposed axiom M is consistent with the axiom of choice, the present argument cannot be reframed as an argument against the axiom of choice.[4] This makes it importantly different from Freiling's chance-based argument against the continuum hypothesis. While Freiling's symmetry axiom $A_{\aleph_0}$ is also compatible with choice, his reasoning in support of that axiom extends to justify stronger axioms conflicting with choice, and arguably even axioms that conflict with ZF — Freiling pointed this out himself. This is generally taken to be the main reason why this argument won so few converts: to accept Freiling's premises one must enter a slippery slope towards denying some highly plausible, fruitful and widely accepted mathematical principles (see Maddy 1988, §II.3.10). By contrast, the argument of this paper is tailor made for M, and cannot without further assumptions be used to argue against ZFC. Another contrast with Freiling is that the present argument has a stronger conclusion, purporting to show that the continuum is greater than various large cardinals.

---

[2] To be precise, Solovay showed ZFC + M is *equiconsistent* with ZFC + "There exists a two-valued measurable cardinal". Most set theorists would agree the latter theory is consistent, so that ZFC + M must be consistent, too.

[3] A *real-valued measurable* cardinal $\kappa$ admits of a total probability measure that is not only countably additive, but also $\kappa$-additive, meaning that if $\lambda < \kappa$, then any union of $\lambda$ sets of measure zero has measure zero itself. While ZFC + M proves the existence of such a cardinal $\kappa$, it does not state that $2^{\aleph_0}$ itself is one — it may be that the continuum has a countably additive measure of the appropriate kind, but no $2^{\aleph_0}$-additive measure. M also does not entail the existence of a *measurable* or *two-valued measurable* cardinal, which is a cardinal $\kappa$ that admits of a $\kappa$-additive, total probability measure that assigns measure 0 or 1 to every subset of $\kappa$. Only real-valued measurable cardinals greater than $2^{\aleph_0}$ are guaranteed to be two-valued measurable. See Solovay 1971; Jech 2002, Ch. 10.

[4] Except, of course, if you are willing to use the continuum hypothesis (CH) as a premise against the axiom of choice (AC): after all, any argument can be refashioned as a *modus tollens* against its premises. In particular, ZF + M + CH ⊢ ¬AC, and Solovay 1970 has a consistency proof for ZF + M + CH. But unlike with CH, there are strong intuitive and abductive reasons to accept AC. So I think only someone skeptical of AC for independent reasons would have any inclination to construe the argument this way.

## I. Improbable Outcomes

Adopting the terminology of Easwaran 2014, let us say an event is *minuscule* if it has a chance that is less than any positive number. Assuming that chances are real numbers, minuscule events are simply zero-chance events. But as discussed, I am happy to leave open the possibility that some minuscule events have infinitesimal chances. For the spinner and the other chancy procedures mentioned, every outcome (that is, every singleton event) is minuscule, and there are a continuum of those outcomes. Thus these procedures are all illustrations of the fact that some outcome spaces can be divided into a continuum of minuscule events. This is a property of the continuum that is not trivially shared by lesser infinite cardinalities. In particular, because of countable additivity, an outcome space can never be divided into $\aleph_0$ minuscule events — more on that in Section III. This paper is about how we can leverage this special fact about the continuum in order to investigate its size.

But first, I should argue that it is indeed a fact, defending premises (A) and (B) above. At the end of the section, I will explain why it looks to be a *special* fact — in that we have no reason to think lesser infinite cardinalities have this property, too. The most perspicuous example for these purposes is the infinite row of flipping coins. Premise (A) is self-evident in that case: if we flip all the coins, the outcome is always some particular infinite sequence of heads and tails. Premise (B) says that no individual sequence has a positive chance. To see this, take an arbitrary sequence. Say it starts HTTHHT… In order to get that sequence, the first coin must land heads, which has chance ½. And the first three coins would have to land HTT, which has chance ⅛. And so on: consider $n$ coins to see that the chance of getting the whole sequence is at most $\frac{1}{2}^n$. So the chance must be less than any positive number. If the coins are not quite fair, one can still argue this way, provided the chances for individual coin results have an upper bound $r < 1$: for any $n$, consider the first $n$ flips to see that the chance of the sequence is less than $r^n$. So it cannot be positive.

To get a continuum-sized outcome space by coin flips alone, you will have to perform infinitely many flips. Given actual physics, these flips are bound to be spread out over an infinite stretch of space or time. As we will see in Section II, special difficulties attach to such procedures. So it is worth our while to defend (A) and (B) in the case of the spinner too. A spinner selects between a continuum of possible angles in a single whirl, so this procedure only requires a finite amount of space and time.

In the spinner variant of the argument, (A) states that the spinner picks out a particular angle over all the others. To show this, we do not need a physically unrealistic idealisation like a completely sharp and perfectly rigid pointer. Perfect precision, in the sense required, is a matter of settling on a clever and flexible way to determine what angle counts as being "picked out". Here is one such way. Let $C$ be the centre of gravity of the spinner. Now take any off-centre physical part of the spinner, and call it the *pointer*. Label the pointer's centre of mass prior to the spin $P$, and its centre of mass after the spin $P'$. Whatever happens, a whirl of the spinner thus determines a particular plane angle $\angle PCP'$ (as seen in diagram 1 above). (If the spinner is assumed to be rigid, this angle $\angle PCP'$ is moreover independent of the choice of pointer, in that any choice of an off-centre pointer yields the same result.) We can visualise the outcome space $\Omega$ of the possible angles as a continuous circle around $C$. The intersection between the ray $CP$ and $\Omega$ is $0°$, and the intersection between $CP'$ and $\Omega$ is the random angle picked out by the spinner. The continuity of $\Omega$ is essential here: to guarantee that the rays $CP$ and $CP'$ both intersect $\Omega$ no matter what happens, we need to assume that the circle is Dedekind complete. And in order for that to be the case, $\Omega$ has to be a continuum.

It may be objected that, due to the vague boundaries of the spinner and its pointer, the recipe I wrote down does not unambiguously designate a particular angle in $\Omega$. But our inability to determinately specify these physical objects does not make the objects themselves indeterminate (Evans 1978). I showed how any pair of a spinner-shaped physical object and a pointer-shaped part of it determines a particular angle on any spin. The argument does not require anyone to specify the members of this pair, merely that a suitable pair exists. If there are many suitable pairs, all the better! It does not matter that they are so similar to one another that it is impossible for us to pick one out determinately.

Premise (B) says that the spinner is *diffuse*, meaning that no particular angle in $\Omega$ has a positive chance of being picked out. There are various ways to show this. One argument starts from the assumption that the spinner has no "supermagnetic" angles, so that no individual angle is infinitely more likely than any other. Since there are uncountably many angles, they cannot all have a positive chance. But then, by our assumption, *no* angle has a positive chance. Here is a different argument. Take an arbitrary point $P$ from $\Omega$. Consider an arc $AB$ including the point $P$. There is some chance $x$ that the spinner lands in $AB$. Since $P$ falls within $AB$, the chance of hitting $P$ is at most $x$. If $x$ is zero or infinitesimal, so is the chance of landing on $P$. Suppose $x$ is positive. If the spinner can land (almost) anywhere in the arc $AB$, it is plausible that a small enough arc $A'B'$ around $P$, with extremities $A'$ and $B'$ closer to $P$, has a chance under $x/2$ of being hit. But then the chance of hitting $P$ is less than $x/2$. By the same token, it is less than $x/4$, less than $x/8$, etcetera. So it is less than any positive number.

All such arguments for premise (B) can be resisted by taking the view that determinism is both metaphysically necessary and incompatible with the existence of non-trivial chances. On that view, there is necessarily only a single angle the spinner could land on, and a single way the coins could land. One might also deny the possibility or coherence of non-trivial chances for other reasons. I grant that, given such views, you should not expect to learn anything about the continuum by thinking about chance. But anti-realism about chance is not a happy position to occupy. The notion of objective chance is very fruitful in science, where probabilistic methods of prediction have met with great empirical success. This fact offers a powerful reason to think the concept of chance has a reflection in reality, or at least in a possible version of it. Moreover, non-trivial chances can arguably arise even in the presence of determinism, as in statistical mechanics — see Loewer 2001 and Glynn 2010.

Are there physically possible, diffuse random processes with outcome spaces smaller than the continuum? It is not clear there are. Certainly, no procedure based on coin flips can give you an outcome space with fewer than $2^{\aleph_0}$ outcomes, all of which are minuscule. For in order to make the outcomes minuscule, you have to perform infinitely many flips (since any finite sequence of heads and tails has a positive chance). But the moment you do infinitely many coin flips, the size of the outcome space balloons to $2^{\aleph_0}$. So there is not going to be any way to divide an outcome space into less than $2^{\aleph_0}$ minuscule outcomes using coin flips, or dice, or anything like that.

What about spinners? Here the matter is more subtle. You might imagine, after all, that there could be an $\aleph_1$-*spinner*: a diffuse spinner that is for some reason physically constrained to land within an $\aleph_1$-sized set of angles. Even if $\aleph_1 < 2^{\aleph_0}$, it seems conceivable that such a spinner could exist. Note that this does not put $\aleph_1$-spinners and $2^{\aleph_0}$-spinners dialectically on a par: assuming diffuse random spinners are possible at all, we have conclusive mathematical reason to think a diffuse spinner with a continuous outcome space is possible. But there is no guarantee for the possibility of $\aleph_1$-spinners. As a matter of mathematical necessity, any spinner must pick out an angle from the continuum $[0, 2\pi)$, since that set is Dedekind complete. But we have no similar, positive reason to support the possibility of a diffuse spinner that can only land within an $\aleph_1$-sized subset of $[0, 2\pi)$.

This paper makes the further point that, moreover, we *do* have good reason to *reject* the possibility of $\aleph_1$-spinners, in light of Banach and Kuratowski's theorem. Thus we are led to accept the possibility of $2^{\aleph_0}$-spinners while rejecting the possibility of $\aleph_1$-spinners, forcing the conclusion that $2^{\aleph_0}$ does not equal $\aleph_1$. In the same way, the other results cited above distinguish the continuum from other infinite cardinalities like $\aleph_2$. It may be worth noting, however, that my premises *are* in fact consistent with the possibility of minuscule-outcome spinners whose outcome spaces have a size $\kappa$ below the continuum. But such a $\kappa$-spinner would only be possible for $\kappa$ that are much greater than $\aleph_1$. To be precise, $\kappa$ would have to be greater than or equal to the least real-valued measurable cardinal (see footnote 2 above). It should not surprise us that set-theoretical results should bear on which $\kappa$-spinners are and aren't possible. After all, metaphysical possibility is constrained by the mathematical facts.

Besides coin-flip based procedures and $\kappa$-spinners, one could try to think up a different way of picking one of $\kappa < 2^{\aleph_0}$ outcomes.[5] However, so far no-one has ever succeeded in actually constructing a convincing method. Particularly in the case $\kappa = \aleph_0$, this is not for lack of trying: John Norton (2018) documents dozens of designs for a so-called *de Finetti lottery* — a fair random method for picking a natural number. But as his discussion shows, none of them ultimately work. Norton also proffers his own suggestion, but was forced to concede in Norton and Pruss 2018 that it, too, is fatally flawed. I reckon this is no coincidence. For if chances are countably additive (Section III), de Finetti lotteries are impossible; and if they are always total (Section II), diffuse procedures with $\aleph_1$ different outcomes are impossible too. But there are diffuse procedures with $2^{\aleph_0}$ outcomes. So $2^{\aleph_0}$ must be distinct from $\aleph_1$.

## II. Totality

Now we arrive at the controversial part: the *Totality* premise (C). Depending on the example, this is the claim that that every spinner event has a chance, the claim that every infinite coin flip event a chance, and so on. Let us say a physical process is *totally chancy* if it has this feature. It is uncontroversial that procedures with finitely many outcomes, like die rolls, are totally chancy. But that assumption is not standardly extended to processes with a continuum of outcomes. In §2.2, I go into the reasons for this, and explain why I think they are bad reasons. But let me start with the positive case for thinking that physical systems, if they are chancy at all, had better be totally chancy.

### 2.1 The Abductive Case for Totality

Imagine you are given a lamp, call it the *Mystery Lamp*, that is operated by a single green button on its base. When you press the green button, sometimes the Mystery Lamp blinks at you after a second. Sometimes it stays dark. The internal operations of the lamp are governed by a random mechanism hidden on the inside: this mechanism is set in motion whenever the button is pressed, and whether or not the lamp blinks depends in some way on the final state of that internal mechanism. While the nature of this random mechanism is entirely unknown to us (hence the lamp's name), we *do* know that it returns to its initial state after each run, so that the results of successive presses are causally independent. We also know that the mechanism is robust and does not change between trials. Without having an opportunity to examine the innards, you are asked to investigate how likely the lamp is to blink at a given press of the button.

---

[5] Here is one design that may have occurred to the reader. Let $\aleph_0 \leq \kappa < 2^{\aleph_0}$. Now pick a $\kappa$-sized set of angles $X$ from $[0, 2\pi)$ and spin a random spinner until eventually we hit an angle in $X$. The chance of ending up on any particular angle $x \in X$ is minuscule, and there are $\kappa$ possible outcomes! Norton (§2.3) considers this design for the case $\kappa = \aleph_0$. But as he observes, it has a fatal bug. We are overlooking an important outcome: it is possible that the procedure fails to terminate, because we only get angles from $[0, 2\pi) \backslash X$ forever. On plausible assumptions, the chance that this happens equals 1, because $[0, 2\pi) \backslash X$ outnumbers $X$ to the tune of $2^{\aleph_0}$ to 1. So on reflection, the outcome space of this procedure is not partitioned into $\kappa$ minuscule cells after all: it has $\kappa$ minuscule cells, and one huge cell of measure 1.

Here is how you would ordinarily go about it. You press the button a number of times, recording each time whether the lamp blinks or not. Suppose you do a hundred trials and record eighty-one blinks. You can conclude with some confidence that the chance of blinking on a given trial is about 80%. That hypothesis explains your findings so far, and it also justifies certain expectations about the future behaviour of the lamp: we have a sophisticated understanding of random behaviour in chancy situations. This is, in the simplest instance, how probabilistic induction works in science. And the core motivation for espousing realism about chance, I take it, is that it promises to explain why this inference method should be so reliable.

Now if the random mechanism inside the mystery lamp is assumed to be totally chancy, we have at least a schematic idea of how that explanation goes. Different hypotheses about the chance of blinking yield different expectations for the blinking frequency in a string of a hundred trials. By observing the actual frequency, we test the full gamut of chance hypotheses, confirming some while disconfirming others. And that helps us determine what behaviour to expect in the future.

But if the mechanism in the lamp were *not* totally chancy, that means we have to contend with an important additional hypothesis: for in that case, it may be that the blinking event has no chance at all. The problem is that if this is so, we have no idea what to expect. The observation of 81 blinks over a hundred trials disconfirms the hypothesis that the chance of blinking is 70%, and very strongly disconfirms the hypothesis that the chance of blinking is 10%. Does it confirm or disconfirm the hypothesis that this event has no particular chance at all? We cannot say, and that is itself highly problematic if we are trying to vindicate ordinary inductive inferences. For if the hypothesis that the event of blinking has no chance fits the data equally well as the hypothesis that the chance of blinking is around 80%, then that substantially undermines the strong justification we thought we had for our inductive predictions about the future behaviour of the lamp.

There are many metaphysical mysteries about the nature of objective chance. Nonetheless, chance is in some sense the most well-understood kind of future contingency there is. When something has a chance, we know what to expect. By contrast, if an event lacks a chance altogether, we have no idea what to expect. Consequently, the existence of chance-free events wreaks havoc on the standard justification of probabilistic inductive reasoning. This happens even if the event under examination in fact has a chance — the mere epistemic possibility, however slight, that it might lack a chance is what gets us into trouble. Our ordinary inductive reasoning about such cases rules out that possibility from the outset; in other words, it presupposes that any event has a chance. That strong, simple assumption provides an elegant justification for ordinary probabilisitic inductive reasoning, which would be undermined by the existence of chance gaps (events that lack a chance).

This issue about induction is not the only problem with the notion of chance gaps: they raise a range of awkward metaphysical and epistemological questions that lack satisfactory answers. Why do some events have a chance while others lack it? How are we to know that a particular event has a chance? If an event has no chance, what should be our credence that it will happen? Can chance gaps arise only in complex, infinite-outcome setups like a spinner? Or could there also be a chance-free coin flips and chance-free dice? If not why not? And if so, how do you make a chance-free die?

These basic questions on chance gaps have received little to no attention in the literature, but Hájek and Smithson 2012 and Hájek, Hawthorne and Isaacs ms. do make progress on some of them. These authors note that, even if an event lacks a precise chance itself, it can still be described by what they call an *indeterminate chance*. Rather than a real number, an indeterminate chance is a closed interval of real numbers between 0 and 1. For a given chance-free event $E$, the lower bound of this interval is its *inner measure* — the supremum of all the chances of events that entail $E$ and have a precise chance. Its upper bound is the outer measure. The view then addresses epistemological issues about chance gaps with a generalised version of Lewis' Principal Principle, making use of the fact that the authors also support the idea of imprecise or "mushy" credences. The generalised principle runs as follows: for any proposition $P$, conditional on admissible evidence entailing that $P$ has an indeterminate chance $[a, b]$, you should have an imprecise credence that $P$ equal to $[a, b]$.

This theory of indeterminate chances makes progress on many questions about chance gaps. But it only brings out the difficulty about probabilistic inductions with greater urgency. In particular, in spite of the generalised Principal Principle, it is still entirely opaque just how our expectations about the future should be affected by the existence of indeterminate chances. One issue here is that imprecise credences are themselves problematic (White 2009, Elga 2010, Schoenfield 2017). But there is another, more pressing difficulty. Recall our investigation of the Mystery Lamp, which found 81 blinks in a string of a hundred trials. What are we to make of this finding in a world with indeterminate chances? The conclusion that there is a determinate chance of around 80% of blinking is presumably no longer warranted. It could just as well be some indeterminate chance. But which one would it be? A [79%, 81%] chance? Or [50%, 90%]? Or [0%, 100%]? How are we even supposed to tell the difference? Is any hypothesis more plausible than the others? And, crucially, what difference does it make for our expectations about the future?

Hájek and friends offer no clear answer. But let us cut them some slack and assume they eventually meet those challenges one way or another. That is, suppose they gave us a way of finding out whether and how indeterminate the chance of a given event is, and a way to extract clear predictions from indeterminate chance hypotheses. Inevitably, those predictions would often differ from the predictions extracted on the basis of our present, cruder inductive methods, which overlook the possibility of indeterminate chances entirely. So the view would motivate a revision of our current probabilistic inductive methods. Will the new methodology produce better science than the present one? Only time will tell, but I hope I will not be thought too much of a spoilsport for sounding a note of scepticism. The basic reason for believing in objective chance is to help justify our actual inductive methods, and to explain their remarkable success. If the metaphysics of chance we end up with recommends a sweeping revision of those very methods, we have clearly gone astray.

The simple, naïve view that all events have a determinate chance recommends itself as a scrutable and conceptually clean metaphysical vision that underwrites probabilisitic reasoning in science as it is actually conducted. Moreover, we should be skeptical that a more sophisticated alternative, which differentiates chancy and chance-free events, or sharp and indeterminate chances, can accomplish this. Standard statistical reasoning is hugely predictive in science, and it builds on the assumption that every event has a chance. This gives us excellent abductive reason to think that assumption is true.

## 2.2  Totality and Symmetry

Why is this view not widely endorsed? The reasons have to do with mathematical practice: standardly, mathematicians will define a measure over a continuum-sized space on a subalgebra of its powerset. For instance, the standard uniform measure on the interval $[0, 2\pi)$ has the class of Lebesgue subsets of $[0, 2\pi)$ as its domain, and assigns no value to non-Lebesgue sets of points. Why do mathematicians ignore non-Lebesgue sets? Is there a theorem showing that spinners have chance gaps? Has it been mathematically proven that there is no such thing as the chance of a non-Lebesgue event? All too often, theoretical treatments of chance suggest the answer is "yes". Van Fraassen 1991 (ch. 3) is a case in point. His discussion is sophisticated and provides a detailed explanation of the relevant mathematical result, Vitali 1905. But then, Van Fraassen simply concludes "Therefore the requirement to have [chance] measures defined everywhere is unacceptable" (p. 55), as if the theorem in question were that total chance functions on the continuum do not exist. But that is not what Vitali showed. And while it is true that this result is the reason that mathematicians restrict their chance measures to Lebesgue sets, it does not give philosophers any clear reason to do the same.

What Vitali *did* show is that, given the axiom of choice, countably additive, total chance functions on a continuum lack a certain kind of symmetry: more specifically, they fail to be translation invariant, or rotation invariant in the case of a circle. Let us say that a spinner is *totally fair* if it has the following property: for every set $S$ of points on the circle, and every rotation $\rho$, the chance of landing on a point in $S$ is equal to the chance of landing on a point in $\rho S$. Vitali's theorem shows that totally fair spinners are not possible. This is an interesting and surprising result. Sure, no-one would argue that any actual, concrete spinning top or fidget spinner was *totally* fair. Realistically, these things always have some small bias or imperfection in some direction or other. So it was clear from the outset that a totally fair spinner is a physically unrealistic idealisation. However, it is striking that this ideal cannot be achieved even in principle.

However, if one restricts the measure Ch to the Lebesgue regions of the circle, then one can at least get a function Ch that satisfies the equation $\text{Ch}(\rho L) = \text{Ch}(L)$ wherever Ch is defined. So by ignoring the non-Lebesgue sets, we gain an important symmetry. In mathematics, symmetries are crucial for reining in complexity: without them, you can't prove theorems, and your models become intractable. So mathematicians ignore non-Lebesgue subsets for the same reason they ignore most things: to achieve the symmetry and simplicity they need in order to do mathematics. This practice implies no positive judgment that these events lack a chance, merely the practical judgement that those chances are best ignored in order to manage the complexity of the situation.

Mathematicians, then, have excellent reason to ignore non-Lebesgue spinner events. But we philosophers do not, at least not if we are interested in the nature of future contingency. For one, it is questionable whether it is even possible draw any sharp general boundary between Lebesgue and non-Lebesgue propositions about the future that corresponds to the mathematical distinction between Lebesgue and non-Lebesgue sets of real numbers. But let's not dwell on that difficulty. Even if we can draw this distinction, surely we should not *ignore* the non-Lebesgue spinner events. The least we can do is ask whether these events do or do not have a chance.

We have already discussed the drawbacks of the negative answer to this question in §2.1. The main problem is that the possibility of chance-free events, however slight, throws the justification of probabilistic induction into serious jeopardy. Meanwhile the positive answer, *Totality*, is relatively simple, perfectly consistent with ZFC, and yields a solid metaphysical underpinning for probabilistic reasoning. The only cost is that *Totality* forces us to give up the intuition that a perfectly symmetric random spinner is possible. Given total chanciness, and our other assumptions, there must be some sets of angles $E$ and rotations $\rho$ such that $\text{Ch}(\rho E) > \text{Ch}(E)$. But how much of a cost is this, really? As noted above, there is little empirical motivation for the contrary intuition: in the real world, spinners never have perfect rotational symmetry.

What is more, related intuitions about rotational symmetry are uncontroversially mistaken. Intuitively, rotations on an ideal spinner should preserve not only the chance of events but also their logical strength, in that it should be impossible to rotate a set of points onto a strict sub- or superset of itself. This is true for finite sets of points, but the principle does not hold in general: most rotations do in fact map some events onto strictly entailed ones. For instance, consider the event

$N = \{\, n \in [0, 2\pi) : n \text{ is a positive integer modulo } 2\pi \,\}$

It is easy to see that if we rotate $N$ counterclockwise by an integer number of radians, the image is a strict subset of $N$. Similarly, we can map $N$ onto a strict superset of itself by rotating the set clockwise.

I do not deny that the view that random spinners are necessarily asymmetric goes against our untutored intuitions. Intuitively, there appears to be no reason why a perfectly symmetric random spinner should be metaphysically impossible. But isn't that appearance just mistaken? In fact there *is* a reason: it is Vitali's theorem. Mathematical facts constrain metaphysical possibility. Surprising mathematical results, like Vitali's, constrain it in surprising ways.

All in all, I think the considerations in favour of total chanciness clearly outweigh those in favour of the possibility of perfect symmetry in the case of a random spinner. In the infinite coin flips case, the balance of reasons seems less decisive in two respects: on the one hand, the abductive argument for total chanciness has less of a foothold, because we have no actual experience of infinite procedures; on the other hand, the intuitive price we pay here is steeper.

For it turns out that, in an infinite array of coin flips, total chanciness rules out the possibility that the flips are both perfectly fair and perfectly independent. Here *perfect fairness* means that the chance of heads exactly equals the chance of tails for every coin, and *perfect independence* means that the chances of two propositions are independent whenever the truth of these propositions supervenes on the outcomes of disjoint sets of coin flips. Independence between individual flips, and between finite sets of flips, *can* be maintained even if the coins are fair — we are not forced to say that the outcome of one flip affects the outcome of another. To be sure, it remains a counterintuitive result. But then again, so it goes: surprising mathematical results constrain metaphysical possibility in surprising ways.

Moreover, a puzzle due to David Builes (2020) gives independent reason to question the tenability of perfect independence in an infinite sequence of fair coin flips.[6] The contrary intuition is, once more, an extrapolation from the symmetries of finite cases — a risky step that has failed us time and again. In any finite array of coin flips, the flips are fair and independent just in case every permutation of the outcome space preserves the chance of an event. But this permutation symmetry fails in the infinite case whether or not we accept *Totality*. Not every symmetry that can be achieved in a finite array of coin flips can be achieved in an infinite array. So on balance, I think there is still a strong case for *Totality* in the infinite coin flip case. Nonetheless, in view of the differences I pointed out, I think it would be reasonable to feel more confident in *Totality* as applied to finite setups like the spinner case. For our ultimate purpose here, that nuance matters little, however. After all, we only really need one version of the *Totality* premise.

---

[6] What is the probability that coin flip #1 in an infinite sequence lands heads, given that only finitely many coins land heads? Builes offers good reasons to think it is less than ½. But this violates perfect independence: whether or not finitely many flips land heads is fully determined by the way all the other coins land.

## 2.3 Hidden Chance Gaps

My main argument for the view that random physical processes are totally chancy has been that scientific reasoning about chance is built on that presupposition. This may not be explicit in the chance functions scientists and mathematicians employ, which ignore some events to allow for symmetry. But it does reveal itself implicitly in our inductive methods, whose justification relies essentially on the assumption that every event has a chance — or so I claim. I now want to consider one strategy for resisting this abductive argument for *Totality*, by maintaining that probabilistic inductive reasoning can in fact be grounded on a weaker assumption, compatible with the existence of chance gaps.

Here is the objection I have in mind. It begins by conjecturing that scientists are only ever talking about events of a certain kind — let's say events that are in some important sense epistemically accessible or *observable*, meaning roughly that we can reliably tell whether or not the event took place.[7] The next step is to explain the success of probabilistic inductive reasoning in science in terms of the chanciness of observable events alone. Consistent with this weaker commitment, the objector is then free to reject total chanciness, and to maintain that chance-free events do exist — it is just that, conveniently, they all happen to be unobservable. In a nutshell: the abductive argument fails to rule out the possibility of *hidden* chance gaps.

Let's make that a bit more concrete. If we assume that the chance-free events in the spinner case are just the non-Lebesgue ones, then the idea that such events are "hidden" in this sort of way has some plausibility. A non-Lebesgue event $N$ contains unnatural, scattered collections of isolated points, which makes it difficult to see how anyone could reliably tell by direct observation whether or not $N$ occurred. Suppose, moreover, that non-Lebesgue events are really entirely unobservable, in that they are causally isolated from all observable events. Then it would be impossible to build a Mystery Lamp which blinks just in case the spinner lands on an angle in $N$. For this reason, the objector continues, the argument from §2.1 does not apply to unobservable, non-Lebesgue events like $N$: our reasoning about the Mystery Lamp does not require us to assume events like $N$ have a chance, because it is physically impossible for such an event to be causally linked to the blinking of the lamp.

---

[7] The informal use of the adjective *observable* in this subsection is unrelated to the quantum mechanical notion of an observable. In quantum mechanics, the exact position of a particle and the exact angle of a spinner pointer technically count as observables, even though there is no practical way to ascertain them.

This is a rather elaborate objection, but it is worth addressing: something in the vicinity of this view comes naturally to those of us who have learnt to live with the idea of chance gaps. One problem with the objection is that it is built on some rather ambitious promissory notes. The hypothesis that scientists only ever reason about observable events stands in clear need of clarification and justification. The same is true for the view that there is a substantial class of unobservable events that exist in perfect causal isolation from those observable events. Another issue is that it is unclear to what extent the objector's explanation for the success of probabilistic induction really holds its own against the simple explanation offered in §2.1. For although this alternative explanation is less committal, it is also a good deal more complex, because it invokes the notion of observability. And as a rule, simple explanations are better than complex ones.

But I will focus here on a different reply, related to this last point. Let's spot this objector the distinction between observable and unobservable events, and let's assume non-Lebesgue event are unobservable. Let's also assume that their alternative account of the success of probabilistic inductive reasoning can be made to work. I contend that, even if we granted all that, the abductive argument against chance gaps maintains its force. For in order to undercut this abductive argument, it is not sufficient for the objector to show they are able to consistently maintain that chance gaps exist. After all, this argument was never shooting for validity. To resist the argument, the objector has to claim that the success of probabilistic inductive reasoning gives us *no good reason* to think all events have a chance. But the objector has not shown this. On the contrary, I would contend that the proposition *all observable events have a chance*, which the objector endorses, makes for a great reason to think unobservable events have a chance too. After all, the class of observable events is very large and varied. So even if we grant the objector their alternative explanation, we still have excellent justification for the thesis that both observable and unobservable events have a chance.

The objector's position is comparable to that of a skeptic who resists the inductive evidence that all emeralds are green, insisting that our observations are fully explained on the weaker hypothesis that *observable* emeralds are green (excuse the artificial example). "As for emeralds that are very far away or too small to be observed with the naked eye," this sceptic proclaims, "I frame no hypotheses! For all I know, they could be pink or navy blue or burgundy with yellow polka dots." Such a level of inductive restraint is clearly excessive. There is no reason to think that emeralds should have a different colour just because they are small or far away. Absent such a reason, we ought to judge it likely that they too are similar to the known sample. We can admit that the conclusion that observable emeralds are green is especially secure, and that it is even better supported than the thesis that all emeralds are green. But the inductive evidence still bears on all emeralds, and the more general conclusion is justified as well.

I think the same thing is true in the present case. The abductive evidence that observable (Lebesgue) events have chances bears on unobservable (non-Lebesgue) events too, because we have no sufficiently compelling a priori reason to expect observable and unobservable events to be different in this respect. This is especially true since chance, whatever it is, seems to a very basic property of events. The contention that there are events that lack this basic property is something that cries out for explanation, particularly when accompanied by the rather suspicious qualification that the events in question happen to be categorically unobservable and causally isolated.

Admittedly, there is a wrinkle in the emerald analogy. While there was no reason at all to think unobservable emeralds should be anything out of the ordinary, the same is arguably not true of non-Lebesgue events. We may have no direct experience of such events, but as Joel Hamkins (2012, 2015) emphasises in his discussion of Freiling, mathematicians do have a lot of distinctively mathematical experience with the sets that represent them. Moreover, they have found that non-Lebesgue sets can behave in strange ways, overturning our ordinary mathematical expectations. Hamkins claims that this is the real source of mathematicians' resistance to Freiling's argument: "We are skeptical of any intuitive or naive use of measure precisely because we know so much now about the various mathematical pitfalls that can arise, about how complicated and badly-behaved functions and sets of reals can be in terms of their measure-theoretic properties." (Hamkins 2012, p. 17)

Hamkins' main point here is about the unreliability of our pre-theoretical intuitions, which does not directly impugn the present line of argument: unlike Freiling, I am building my case on broadly abductive, empirical considerations and not on intuitions. But I think that the considerations Hamkins adduces are relevant to the abductive argument against chance gaps as well. Thanks to Solovay's consistency proof, cited above, we can rest assured that none of the odd formal properties to which Hamkins alludes *logically* entail the existence of chance-free events. Still, one might worry that non-Lebesgue events might somehow lack a chance in virtue of their strange formal properties. The lack of homogeneity between observed and unobservable events is reason to be cautious in generalising to the class of all events. It is, after all, uncontroversial that inductive generalisations over heterogeneous domains are less secure. For instance, an inductive generalisation over the class of all mammals is, all things being equal, riskier than a generalisation over the class of all magpies.

Now I concede that the known heterogeneity of the class of events weakens the strength of the abductive case for the totality premise (C). But I think it would be a serious overreaction to dismiss the argument on those grounds alone. Sure, the fact that some non-Lebesgue events are strange in other respects lends *some* credence to the idea that they may also instantiate the strange property of lacking a chance. But that is just to say that the existence of chance-free events would be *even more* mysterious if non-Lebesgue events were perfectly normal and well-behaved in all other respects. It does not erase the fact that chance-freeness is another very peculiar property, whose instantiation should puzzle and surprise us. Moreover, to reiterate, Solovay's consistency result tells us that this really would be an additional surprise, logically independent of the curious properties ZFC establishes.

It is instructive to view the present conundrum in the context of the debate around the use of inductive justification in other domains of mathematics, like number theory. As Frege once observed, the domain of natural numbers has "none of that uniformity, which in other fields can give [inductive reasoning] a high degree of reliability. … Position in the number series is not a matter of indifference like position in space … each [number] is formed in its own special way and has its own unique peculiarities." (Frege 1953 [1883], §10)  It is tempting, on these grounds, to endorse a radically sceptical view of the use of inductive arguments to justify general conclusions in number theory — Alan Baker (2007) explores this position.

But as Baker is quickly forced to acknowledge, it is difficult to bring that position into harmony with the pivotal role that broadly inductive reasoning plays in mathematical practice.[8] Number theorists make reliable judgments about the plausibility of unproven conjectures all the time. Often, those judgments are informed by simple "experiments" like calculating a few instances. In recent decades, these inductive aspects of mathematical inquiry have even become institutionalised with the rise of *experimental mathematics*, a now flourishing branch of mathematical inquiry with its own institutes, journals and conferences (Borwein and Bailey 2008, Baker 2008).

As James Franklin (1987) asks, if those experimental, non-deductive methods truly are unreliable, how are we to account for some maths professors' ability to reliably identify tractable, plausible and unproven conjectures for their PhD students to work on? Timothy Gowers has another nice illustration of the point, writing that "Probably all serious mathematicians believe that, in the long run, each of the digits 0 to 9 occurs in the decimal expansion of π about 10 per cent of the time," even though "Nobody has found a proof, or even an informal argument, that [π is not biased], and nobody

---

[8] By the end of the paper, Baker has softened his position enough to endorse, in §6, an inductive argument for the truth of Goldbach's conjecture. Officially, Baker hangs on to his sceptical view about induction in the realm of numbers, but only by espousing an unusually narrow conception of inductive reasoning.

expects to." Instead, this belief is at least in part based on the fact that the millions of digits of $\pi$ that have been calculated show no bias, and pass various statistical tests for randomness. "Nobody has ruled out the possibility that every digit after the $10^{100}$th is either a 7 or an 8. And yet it is clearly ludicrous to suppose that that might be the case." (Gowers 2007, p. 34)

In line with Baker 2007, one could criticise Gowers' statistical argument about the decimal expansion of $\pi$ by complaining that the finite sample on which it is based is dwarfed by the infinite totality of $\pi$'s decimals. One could also point out that the sample is biased, in that only the initial segment of $\pi$ was considered. But that only shows that in spite of such limitations, an inductive argument can still be very strong. It is not to deny they are limitations: if God handed us statistical analyses of an infinite sample of decimals of $\pi$, with instances taken from evenly spaced stretches along the entire length of the decimal expansion, our inductive evidence would be stronger still. But the theoretical possibility of even better inductive evidence does not imply that the evidence we have is bad. (For further criticism of Baker's position, see Waxman 2017, §3.4.1.)

In a real-world experimental situation there are always *some* properties that divide the unobserved instances from the observed ones, and we can almost always improve our case with a larger, more varied sample. But the success of experimental reasoning inside and outside of mathematics teaches us that inductive generalisations over heterogenous domains can still be highly reliable. So the heterogeneity charge suggested by Hamkins' remarks does not,  after all, pose a very serious threat to the abductive case against chance-free events. I conclude that our empirical evidence supports the strong conclusion that all events have a chance, observable or not. Of course the evidence is not entirely conclusive — inductive evidence never is. All I seek to establish here is that we have good, if defeasible, empirical reasons to think the continuum hypothesis is false.

## III. Countable Additivity

That leaves premise (D), the assumption that chances form a countably additive measure. What may be challenged here is the view that chances are *countably* and not just finitely additive. The assumption of countable additivity is arguably justified by the demand to sustain standard scientific practice. Scientists routinely appeal to the theory of Lebesgue integration and to laws of large numbers, which rely on countable additivity. Still, many prominent theorists have disputed the countable additivity of probabilities, and it would be remiss of me to pass over the controversy in silence (e.g. De Finetti 1990, Savage 1972, Arntzenius et al. 2004, Schurz and Leitgeb 2008, Dorr 2010, Zhang ms).

Given the tolerant attitude I have taken to the notion of infinitesimal chances, I should clarify that the thesis to be defended in this section is the view that the *real parts* of the chances are countably additive.

Advocates of infinitesimal chances sometimes employ a more permissive sense of countable additivity, which allows countably many infinitesimals to "add up" to a positive value (for instance Benci et al. 2013). Countable additivity in this permissive sense would be insufficient for the main argument to go through. So I employ a restrictive use of the term, on which such distributions still count as violations of countable additivity.

The skeptics of countable additivity I just mentioned have for the most part been animated primarily by arguments that subjective probabilities or *credences* are countably non-additive. But the countable additivity of chance is quite a different matter.[9] In particular, by far the strongest reason for thinking that rational credences are not countably additive is that, in certain situations, it seems rational to divide your credences equally between a denumerable infinity of mutually incompatible possibilities. Then your credence in each possibility is 0 or infinitesimal, and yet your credence that one of them is true would be 1. To pose an analogous challenge to the countable additivity of chances, one would need a random process with denumerably many incompatible outcomes that are all equally likely: a de Finetti lottery. But as I already discussed at the end of Section I, no-one has ever identified a plausible real or possible candidate for such a procedure (Norton 2018).

Secondly, there is a strong conceptual reason to think this is no coincidence, and that de Finetti lotteries and other countably non-additive chance distributions are in fact impossible. This is related to the fact that countably non-additive chance distributions uniquely have the property that there is a partition of the outcome space relative to which they make *every* outcome strictly less likely than some other probability function does.[10] Now *that*, it seems, is impossible, or at the very least it would be quite odd (as noted by Dubins 1975, Kadane et al. 1996, Easwaran 2013, Pruss 2013). For suppose $Ch(E) < Ch^*(E)$ for every event $E$ in some partition of the outcome space. Then if you tested the hypothesis that the chance distribution is $Ch$ against the hypothesis that it is $Ch^*$, the latter will be

---

[9] The Principal Principle (Lewis 1980) demands that rational agents apportion their credences to the chances. But it does not follow that the countable non-additivity of rational credences implies the countable non-additivity of chances. If a Principal Principle compliant agent divides their credences between countably additive chance hypotheses, their credences about the outcomes are not guaranteed to be countably additive unless their credences about those hypotheses are countably additive in the first place.

[10] The general result is the following:

> Let Pr be a finitely additive probability. Then there is a probability $Pr^*$ and a partition $\{ E_i \ : \ i \in I \}$
> of the event space such that for all $i \in I$, $Pr(E_i) < Pr^*(E_i)$ if and only if Pr is not countably additive.

*Proof.* Suppose Pr, $Pr^*$ and $\{ E_i \ : \ i \in I \}$ are as specified. Then $\Sigma_{i \in I} Pr(E_i) < \Sigma_{i \in I} Pr^*(E_i) \le 1$. So Pr is not #$I$-additive, where #$I$ is the cardinality of the partition. To show #$I = \aleph_0$. If #$I$ were finite, then by finite additivity $\Sigma Pr(E_i) = \Sigma Pr^*(E_i) = 1$. So #$I$ is infinite. Next, note that $\delta_i := Pr^*(E_i) - Pr(E_i) > 0$ for all $i \in I$. And $\Sigma_{i \in I} \delta_i$ is finite because $\Sigma \delta_i \le \Sigma Pr^*(E_i) \le 1$. Hence #$I$ cannot be uncountable: uncountably many positive values do not add up to a finite value. So #$I = \aleph_0$. Conversely, if Pr is countably non-additive, there is a partition $\{E_1, E_2, E_3 \ldots\}$ such that $\Sigma_n Pr(E_n) = 1 - \varepsilon$ with $\varepsilon > 0$. Now let $Pr^*(E_n) = Pr(E_n) + \varepsilon / 2^n$, and note that for all $n$, $Pr^*(E_n) > Pr(E_n)$. ∎

confirmed over the former *no matter what happens*. Test it again, and Ch\* is confirmed still more. Why would anyone claim that the chance distribution is Ch, if the distribution Ch\* is a better fit for any observations you could possibly make?

This can be brought out more clearly with a concrete example. Suppose a friendly and reputable shaman tells you that he has a meditative ritual for randomly selecting positive integers, which has a countably non-additive chance distribution. In particular, he claims the chance of drawing the number 1 is ¼; the chance of drawing a 2 is ⅛; and in general the chance of drawing the number $n$ is $\frac{1}{2}^{n+1}$. These chances are countably non-additive, because they only add up to ½. As an alternative to the shaman's story, consider this additive chance hypothesis: the chance of drawing a 1 is ½; the chance of drawing a 2 is ¼, and in general the chance of drawing the number $n$ is $\frac{1}{2}^{n}$. These chances are countably additive since they add up to 1. The latter chance distribution definitely is possible, because it is the chance distribution associated with the following random method for picking a positive integer: toss a fair coin until it lands heads, and let the number of tosses it takes be the number selected. Skeptical yet open-minded, you ask the shaman if he would be willing to perform his ritual for you a few times, and he agrees. He closes his eyes and begins to hum. Will the shaman's non-additive divinations be vindicated, or will the additive chance distribution win out?

No clairvoyance is needed to see what the result will be. The first time the shaman opens his eyes, he says "one". Well, that result is twice as likely given the countably additive chances than it would be on the shaman's chances. The next time, he says "three." On the additive hypothesis, there is a one in eight chance to get that number; with shamanic chances, only one in sixteen. Another clear win for the additive chances. And so on. In fact, why even perform the tests? *Any* outcome is twice as likely with the additive chances than with shamanic chances. The shaman's story loses out *a priori*, because every possible result disconfirms it 2:1. Furthermore, any two results disconfirm it 4:1, any three results 8:1, and so on. As I demonstrate in an appendix, every countably non-additive chance hypothesis displays this self-defeating predetermination: our credences in them are rationally constrained to decrease exponentially as we learn the results of more runs of the chancy procedure in question, no matter what those results might be.

The point is not just that it is irrational to believe the shaman. Given the close connection between chance and inductive reasoning, the fact that we cannot say what results would confirm his story indicates that we lack a clear conception of what it would even be for that story to be true. Thus the fact that our ordinary reasoning about chances leaves no room for countably non-additive chance distributions is a good reason to think such distributions are in fact impossible.

⁎⁎

Here is a neat way to sum up the main findings of this paper. Our ordinary probabilistic inductive practices assume that chance distributions are total, and also that they are countably additive. Once we take those assumptions on board, insights in set theory yield all sorts of interesting restrictions on the kinds of chancy processes that are in principle possible. Vitali's result tells us that there can be no rotationally symmetric random method for picking out a point on a circle. The Banach-Tarski paradox reveals that there can be no spherically symmetric random process for picking out an arbitrary point in a continuous sphere. And Ulam's result shows that it is not possible to select one out of $\aleph_1$ possibilities in a way that gives no outcome a positive chance. Since we know independently that it *is* in fact possible to select an arbitrary point from the continuum in such a way (spinners, random darts, infinite sequences of coin tosses), it follows that the continuum cannot be equal to $\aleph_1$.

### References

Arntzenius, Frank, Adam Elga and John Hawthorne, 2004, "Bayesianism, Infinite Decisions and Binding." *Mind* 113(450): 251-283.

Baker, Alan,     2007, "Is there a problem of induction in mathematics?" In *Mathematical Knowledge*, eds. Mary Leng, Alexander Paseau and Michael Potter, 59-73. Oxford: Oxford University Press.
2008, "Experimental Mathematics." *Erkenntnis* 68: 331–344.

Banach, Stefan and Casimir Kuratowski, 1929, "Sur une généralisation du problème de la mesure." *Fundamenta Mathematicæ* 14: 127-131.

Benci, Vieri, Leon Horsten and Sylvia Wenmackers, 2013, "Non-Archimedean Probability." *Milan Journal of Mathematics* 81(1): 121-151

Borwein, Jonathan and David Bailey, 2008, *Mathematics by Experiment: Plausible Reasoning in the 21st Century*. Second Edition. Boca Raton, FL: Taylor and Francis.

Builes, David,     2020, "A Paradox of Evidential Equivalence." *Mind* 129(513): 113-27.

Cohen, Paul,     1963, "The Independence of the Continuum Hypothesis." *Proc. of the National Academy of Sciences of the USA*, 50(6): 1143–1148.

Dorr, Cian,     2010, "The Eternal Coin: a puzzle about self-locating conditional credence" *Philosophical Perspectives* 24: 189–205.

Dubins, Lester E.,     1975, "Finitely additive conditional probabilities." *Annals of Probability* 3: 89-99.

Easwaran, Kenny,     2013, "Why Countable Additivity?" *Thought* 2: 53-61.
2014, "Regularity and Hyperreal Credences." *Philosophical Review* 123(1): 1-41.

Elga, Adam,     2010, "Subjective Probabilities should be Sharp." *Philosophers' Imprint* 10.5, 1-11.

Evans, Gareth,     1978, "Can there be Vague Objects?" *Analysis* 38(4): 208.

Finetti, Bruno de,     1990, *Theory of Probability*. Two volumes, New York: Wiley.

Franklin, James,     1987, "Non-Deductive Logic in Mathematics." *Brit. J. Phil. Sci.* 38(1): 1-18.

Frege, Gottlob,     1953 [1884], *Foundations of Arithmetic*. Translation J.L. Austin, Second Revised Edition. New York: Harper & Brothers.

Freiling, Chris,     1986, "Axioms of symmetry: throwing darts at the real number line.C" *J. of Symbolic Logic*.

Glynn, Luke,     2010, "Deterministic Chance." *British Journal for the Philosophy of Science* 61(1): 51-80.

Gödel, Kurt,     1938, "The Consistency of the Axiom of Choice and of the Generalized Continuum-Hypothesis." *Proc. of the National Academy of Sciences of the USA*, 24(12): 556–557.

Gowers, W.T.,     2007, "Mathematics, Memory and Mental Arithmetic." In: *Mathematical Knowledge*, eds. Mary Leng, Alexander Paseau and Michael Potter, 33-58. Oxford: Oxford University Press.

Hájek, Alan, John Hawthorne and Yoaav Isaacs, "A Measure-Theoretic Vindication of Imprecise Credence." Unpublished manuscript.

Hájek, Alan and Michael Smithson, 2012, "Rationality and indeterminate probabilities." *Synthese* 187(1): 33-48.

Hamkins, Joel David,     2012, "The Set-Theoretic Multiverse." *Review of Symbolic Logic* 5(3): 416-449.
2015, "Is the Dream Solution of the Continuum Hypothesis Attainable?" *Notre Dame Journal of Formal Logic* 56(1): 135-145.

Hofweber, Thomas and Ralf Schindler, 2016, "Hyperreal-Valued Probability Measures Approximating a Real-Valued Measure." *Notre Dame Journal of Formal Logic* 57(3): 369-374.

Honzik, Radek,    2017, "Large cardinals and the Continuum Hypothesis." In: *The Hyperuniverse Project and Maximality*, 205-226. Eds. Antos, Friedman, Honzik and Ternullo, Cham: Birkhäuser.

Howson, Colin,    2014, "Finite additivity, another lottery paradox and conditionalisation." *Synthese* 191: 989–1012.

Jech, Thomas,    2002, *Set Theory*. Third edition. Berlin: Springer.

Kadane, Joseph, Mark Schervish and Teddy Seidenfeld, 1996, "Reasoning to a Foregone Conclusion." *Journal of the American Statistical Association*, 91(435): 1228-1235.

Levy, Azriel and Robert M. Solovay, 1967, "Measurable cardinals and the continuum hypothesis." *Israel Journal of Mathematics* 5: 234-248.

Lewis, David,    1980, "A Subjectivist's Guide to Objective Chance." In: *Philosophy of Probability: Contemporary Readings*, 83-132. Ed. Richard Jeffrey, University of California Press.

Loewer, Barry,    2001, "Determinism and Chance." *Studies in the History and Philosophy of Modern Physics* 32(4): 609–620.

Maddy, Penelope,    1988, "Believing the Axioms, part I." *Journal of Symbolic Logic* 53(2): 481-511.

Norton, John D.,    2018, "How to build an infinite lottery machine." *European Journal for Philosophy of Science* 8(1): 71–95.

Norton, John D. and Alexander R. Pruss, 2018, "Correction to John D. Norton 'How to build an infinite lottery machine'." *European Journal for Philosophy of Science* 8: 71–95.

Pruss, Alexander R.,    2013, "Infinitesimals Are Too Small for Countably Infinite Fair Lotteries." *Synthese*.

Rittberg, Colin J.,    2015, "How Woodin changed his mind: new thoughts on the Continuum Hypothesis." *Archive for History of Exact Sciences* 69: 125-151.

Savage, Leonard J.,    1972, *The Foundations of Statistics*. Second Revised Edition. New York: Dover Publications.

Schoenfield, Miriam,    2017, "The Accuracy and Rationality of Imprecise Credences." *Noûs* 51(4): 667-685.

Schurz, Gerhard and Hannes Leitgeb, 2008 "Finitistic or Frequentistic Approximation of Probability Measures with or without σ-Additivity." *Studia Logica* 89(2): 257-283.

Solovay, Robert M.,    1970, "A model of set-theory in which every set of reals is Lebesgue measurable." *Annals of Mathematics* 92: 1-56.

1971, "Real-Valued Measurable Cardinals". In *Axiomatic Set Theory*, ed. Dana Scott, Proceedings of Symposia in Pure Mathematics 13(1): 397-428.

Ulam, Stanislaw,    1930, "Zur Masstheorie in der allgemeinen Mengenlehre". *Fundamenta Mathematicæ* 16: 140-150.

Van Fraassen, Bas,    1991, *Quantum Mechanics: An Empiricist View*. New York: Oxford University Press.

Vitali, Giuseppe    1905, *Sul problema della misura dei gruppi di punti di una retta*. Bologna: Tip. Gamberini e Parmeggiani. Repr. in Vitali, *Opere sull'analisi reale e complessa*. Florence: Cremonese.

Waxman, Daniel,    2017, *Freedom, Truth, and Consistency*. PhD Dissertation, New York University.

White, Roger,    2009, "Evidential Symmetry and Mushy Credence." In *Oxford Studies in Epistemology*, eds. Tamar Szabo Gendler and John Hawthorne, 161-186. New York: Oxford University Press.

Williamson, Timothy,    2007, "How probable is an infinite sequence of heads?" *Analysis* 67(3): 173-80.

Woodin, Hugh,    2001, "The Continuum Hypothesis, Parts I and II." *Notices of the American Mathematical Society* 48(6-7): 567–576, 681–690.

2010, "Strong axioms of infinity and the search for V." *Proceedings of the international congress of mathematicians*, Hyderabad, India.

2019, *Cantor's Continuum Hypothesis*. Presentation at the Linde Hall Inaugural Maths Symposium at Caltech, delivered February 23rd 2019, youtu.be/hk0iPpeFNLU.

Zhang, Snow,    "A New Argument against Countable Additivity," manuscript.

## Appendix A. Banach and Kuratowski's Theorem

This appendix exhibits a proof of the core mathematical result driving this paper, namely that ZCF + M ⊢ $2^{\aleph_0} > \aleph_1$, where M is the proposed new axiom:

M:  Continuum-sized sets $\Omega$ admit of a total, countably additive measure Ch such that $\mathrm{Ch}(\{x\}) = 0$ for every $x \in \Omega$.

This proof is a version of the original due to Banach and Kuratowski (1929). Adapted to the present context of chance measures, it takes on an especially intuitive form. As mentioned in the text, this result was strengthened by Ulam (1930) and again by Solovay (1971) to show that in ZFC + M, the continuum has to be *much* larger than $\aleph_1$. Ulam's results are covered in Jech 2002 (ch. 10, 131-3).

**Definition**. Let $\mathbf{N}^+$ = { 1, 2, … } be the set of all positive integers. Then if $f$ and $g$ are functions from $\mathbf{N}^+$ to $\mathbf{N}^+$, let us say $f$ is ***smaller than*** $g$, written $f < g$, just in case $f(n) < g(n)$ for all but finitely many natural numbers $n \in \mathbf{N}^+$. (It is easily checked that < is a strict partial order.)

Recall that, in addition to being the first uncountable cardinal, $\aleph_1$ is also the set of all the countable ordinals. Using this fact, we obtain the following lemma.

**Lemma**. If $2^{\aleph_0} = \aleph_1$, then there is a continuum-sized set $\Omega$ of functions from $\mathbf{N}^+$ to $\mathbf{N}^+$ such that for any function $f$ from $\mathbf{N}^+$ to $\mathbf{N}^+$, all but countably many members of $\Omega$ are bigger than $f$.

*Proof*. The cardinality of the set of all functions from $\mathbf{N}^+$ to $\mathbf{N}^+$ is $2^{\aleph_0}$, so if $2^{\aleph_0} = \aleph_1$, we can let { $f_\alpha : \alpha \in \aleph_1$ } be an enumeration of all those functions by the countable ordinals. Now we can construct the desired set $\Omega$ = { $g_\alpha : \alpha \in \aleph_1$ } using the following recursion:

$g_0 \quad : n \longmapsto 1$

$g_{\alpha+1} : n \longmapsto f_\alpha(n) + g_\alpha(n)$

$g_\lambda \quad : n \longmapsto \sum_{k=0}^{n} g_{\beta_k}(n) \qquad$ where { $\beta_k : k \in \mathbf{N}^+$ } enumerates the ordinals below $\lambda$

Since these functions have positive values, it follows from the second clause that for all $\alpha$, $g_\alpha < g_{\alpha+1}$. At limit ordinals $\lambda$, use the third clause to see that $g_{\beta_m} < g_\lambda$ for any ordinal $\beta_m < \lambda$, because for all $n > m$, $g_{\beta_m}(n) < \sum_{k=0}^{n} g_{\beta_k}(n) = g_\lambda(n)$. So { $g_\alpha : \alpha \in \aleph_1$ } is a strictly increasing sequence. Now from the second clause we also know that for any function $f_\alpha$, $f_\alpha < g_{\alpha+1}$. Combining these two facts, we can conclude that $f_\alpha < g_\beta$ for all $\beta > \alpha$, which is to say for all but countably many $g_\beta \in \Omega$. ∎

Suppose there were a set $\Omega$ of this kind, and suppose we represented each of its members with a point on a continuous dartboard at which we are about to throw a random dart. Then that dartboard would have the feature that for any function $f : \mathbf{N}^+ \to \mathbf{N}^+$ whatsoever, the chance that the function $g$ that is hit by the dart is bigger than *f* would be equal to 1, and the chance that *g* is smaller than *f* would be equal to 0. To establish Banach and Kuratowski's result, we show that, assuming M is true, no such dartboard is possible: given any dartboard $\Omega$, we can identify a function $f_\Omega$ that increases so fast that the chance of hitting a smaller function is greater than ½.

**Theorem (Banach and Kuratowski).** ZFC + M $\vdash$ $2^{\aleph_0} > \aleph_1$

*Proof.* Let $\Omega$ be any continuum-sized set of functions from $\mathbf{N}^+$ to $\mathbf{N}^+$. Given M, we may let Ch be a countably additive, total probability measure on $\Omega$ such that $\text{Ch}(X) = 0$ for every countable $X \subseteq \Omega$. We identify a function $f_\Omega$ such that the set of functions smaller than $f_\Omega$ has positive measure, and is therefore uncountable. Since such a function $f_\Omega$ can be constructed for any $\Omega$, it follows from the lemma that $2^{\aleph_0} > \aleph_1$.

Suppose we are about to select a random function *g* from $\Omega$. By countable additivity, we have $\sum_{k=1}^{\infty} \text{Ch}(g(1) = k) = \text{Ch}(\Omega) = 1$, whence there must be a (smallest) integer *K* such that

$$\text{Ch}(g(1) \le K) \;=\; \text{Ch}(g(1) = 1) + \text{Ch}(g(1) = 2) + \ldots + \text{Ch}(g(1) = K) \;>\; \tfrac{3}{4}$$

Now we set $f_\Omega(1) = K + 1$ so that $\text{Ch}(g(1) \ge f_\Omega(1)) < \tfrac{1}{4}$. In a similar way, we can set $f_\Omega(2)$ at a value such that $\text{Ch}(g(2) \ge f_\Omega(2)) < \tfrac{1}{8}$. And in general, for any *n*, we can find an integer $f_\Omega(n)$ such that $\text{Ch}(g(n) \ge f_\Omega(n)) < \tfrac{1}{2^{n+1}}$, as illustrated below. Now note that

$$\begin{aligned}
\text{Ch}(g < f_\Omega) \;&\ge\; \text{Ch}(\textstyle\bigwedge_{n=1}^{\infty} g(n) < f_\Omega(n)) \\
&\ge\; 1 - \textstyle\sum_{n=1}^{\infty} \text{Ch}(g(n) \ge f_\Omega(n)) \\
&>\; 1 - \textstyle\sum_{n=1}^{\infty} \tfrac{1}{2^{n+1}} \\
&=\; \tfrac{1}{2}. \;\blacksquare
\end{aligned}$$

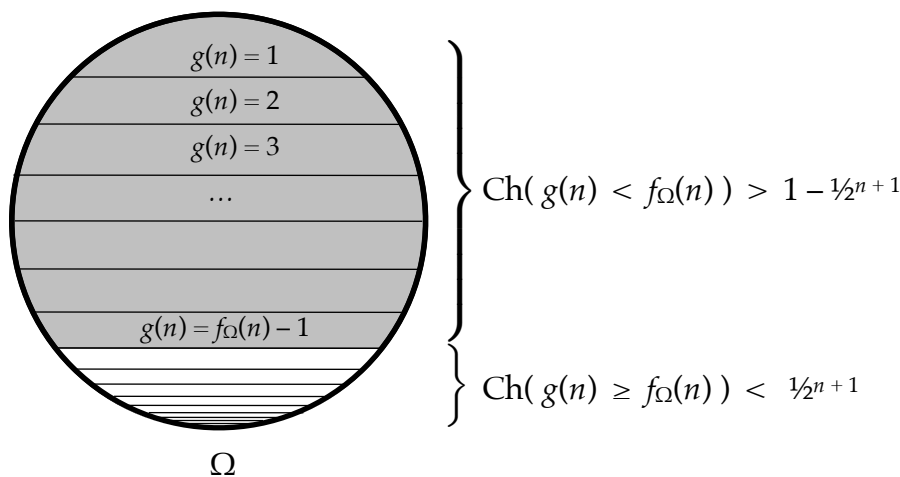

*Diagram 2. For any choice of n, the set of functions $\Omega$ (the "dartboard") is partitioned into infinitely many cells $g(n) = 1$, $g(n) = 2$, … We can always set $f_\Omega(n)$ large enough to ensure that the shaded region $g(n) < f_\Omega(n)$ has measure greater than $1 - \tfrac{1}{2^{n+1}}$, and its complement $g(n) \ge f_\Omega(n)$ has measure less than $\tfrac{1}{2^{n+1}}$.*

## Appendix B. Countably Non-Additive Chance Hypotheses

We compare an arbitrary countably non-additive chance hypothesis $H$ to a suitable alternative hypothesis $H^*$. We then show that $H^*$ is strictly better supported than $H$ no matter what the empirical results are, which justifies the standard scientific practice of ruling out countably non-additive chance hypotheses a priori. Let $H$ be the following claim:

> For any outcome $E_n$ in { $E_1$, $E_2$, … }, the chance it will take place is equal to $\mathrm{Ch}(E_n)$.      (H)

where $\{E_1, E_2, …\}$ is a countable partition of the outcome space, Ch is a finitely additive probability such that $\sum_{n=1}^{\infty} \mathrm{Ch}(E_n) = 1 - \varepsilon$, and $\varepsilon > 0$. Then consider this countably additive variant of $H$:

> For any outcome $E_n$ in { $E_1$, $E_2$, … }, the chance that it will take place is equal to $\mathrm{Ch}^*(E_n) = (1 + \varepsilon) \cdot \mathrm{Ch}(E_n) + \frac{\varepsilon^2}{2^n}$.      (H*)

Note that for every $n$, $\mathrm{Ch}^*(E_n) > (1 + \varepsilon) \cdot \mathrm{Ch}(E_n) \geq \mathrm{Ch}(E_n)$, so that $H^*$ will be confirmed over $H$ by any outcome $E_n$. And note that the chance hypothesis $H^*$ is countably additive over { $E_1$, $E_2$, … }:

$$
\begin{aligned}
\sum_{n=1}^{\infty} \mathrm{Ch}^*(E_n) &= (1 + \varepsilon) \cdot \sum_{n=1}^{\infty} \mathrm{Ch}(E_n) + \varepsilon^2 \cdot \sum_{n=1}^{\infty} \tfrac{1}{2}^n \\
&= (1 + \varepsilon) \cdot (1 - \varepsilon) + \varepsilon^2 \cdot 1 \\
&= (1 - \varepsilon^2) + \varepsilon^2 = 1
\end{aligned}
$$

To compare how $H$ and $H^*$ are confirmed by the evidence, we will make three assumptions:

i) *Open-mindedness*: The agent's prior credence in $H^*$ is greater than zero.

ii) *Bayesianism*: The agent's credences at any time form a finitely additive probability Cr. Suppose the agent has prior credences Cr and receives evidence $E$ such that $\mathrm{Cr}(E) > 0$. Then their posterior credence in any given proposition $X$ is $\mathrm{Cr}(X \mid E) = \frac{\mathrm{Cr}(X \wedge E)}{\mathrm{Cr}(E)}$.

iii) *Principal Principle*: Suppose a draw of the experiment is about to be conducted, leading to one of the outcomes $E_1$, $E_2$ … Then provided $\mathrm{Cr}(H)$, $\mathrm{Cr}(H^*) > 0$, the agent has conditional credences $\mathrm{Cr}(E_k \mid H) = \mathrm{Ch}(E_k)$ and $\mathrm{Cr}(E_k \mid H^*) = \mathrm{Ch}^*(E_k)$ for any $k$.

Claims (ii-iii) are widely accepted rationality constraints (though Colin Howson (2014, §8) avoids the present argument by restricting conditionalisation). Claim (i) is not mandated by rationality, but by a more practical constraint: we wish to compare the degree of confirmation of $H$ and $H^*$, and Bayesian confirmation theory only allows us to do that for hypotheses with non-zero priors.

**Lemma**. Let $H$ and $H^*$ be as above, and suppose the agent satisfies (i-iii). Let $p$ and $p^*$ be the agent's prior credence in $H$ and $H^*$ respectively, conditioned on $(H \vee H^*)$. Suppose a run of the experiment is done, and the agent learns which outcome $E_1$, $E_2$, … has resulted. Then no matter what that outcome is, their posterior credence in $H$ conditioned on $(H \vee H^*)$ is at most $\lambda \cdot p$, where $\lambda = (1 + p^* \cdot \varepsilon)^{-1} < 1$.

*Proof*. Let Cr be the agent's prior credence function conditioned on $H \vee H^*$. By *Bayesianism*, Cr is a finitely additive probability with $\mathrm{Cr}(H \vee H^*) = 1$. By *Open-mindedness*, we have that $\text{Cr}(H^*) = p^* = 1 - p > 0$. Suppose the result of the draw is $E_k$. By *Bayesianism*, the posterior credence in $H$ conditional on $(H \vee H^*)$ will equal $\text{Cr}(H | E_k)$. If $\text{Ch}(E_k) = 0$, it is easily seen that $\text{Cr}(H | E_k) = 0 \le \lambda \cdot p$. If $\text{Ch}(E_k) > 0$, then using the *Principal Principle*,

$$
\begin{aligned}
\text{Cr}(H | E_k) &= \text{Cr}(H) \cdot \frac{\text{Cr}(E_k | H)}{\text{Cr}(E_k)} = \text{Cr}(H) \cdot \frac{\text{Cr}(E_k | H)}{\text{Cr}(H) \cdot \text{Cr}(E_k | H) + \text{Cr}(H^*) \cdot \text{Cr}(E_k | H^*)} \\[2mm]
&= p \cdot \frac{\text{Ch}(E_k)}{p \cdot \text{Ch}(E_k) + p^* \cdot \text{Ch}^*(E_k)} = p \cdot \frac{\text{Ch}(E_k)}{p \cdot \text{Ch}(E_k) + p^* \cdot (1 + \varepsilon) \cdot \text{Ch}(E_k) + p^* \cdot \frac{\varepsilon^2}{2^k}} \\[2mm]
&< p \cdot \frac{\text{Ch}(E_k)}{p \cdot \text{Ch}(E_k) + p^* \cdot (1 + \varepsilon) \cdot \text{Ch}(E_k)} = p \cdot \frac{1}{p + p^* \cdot (1 + \varepsilon)} = p \cdot \frac{1}{1 + p^* \cdot \varepsilon} \\[2mm]
&= p \cdot \lambda \quad \blacksquare
\end{aligned}
$$

**Result 1**. As more draws are performed and an agent satisfying conditions (i-iii) learns the outcomes, that agent's credence in $H$ will exponentially tend to 0, no matter what the outcomes of the draws are.

> *Proof.* Let $p_n$ and $p_n^*$ be the agent's credence in $H$ and $H^*$ respectively, conditioned on $(H \vee H^*)$, and after $n$ runs of the experiment. From the lemma, $p_{k+1} \le \lambda_k \cdot p_k$ for every $k$, where $\lambda_k = (1 + p_k^* \cdot \varepsilon)^{-1}$. Since the sequence $p_k^*$ is monotonically increasing, $\lambda_k$ is decreasing, so that $p_{k+1} \le \lambda_k \cdot p_k \le \lambda_0 \cdot p_k$. By induction, it follows $p_n \le \lambda_0^n \cdot p_0$. The agent's unconditional credence in $H$ is less than that, and is therefore also bounded above by $\lambda_0^n \cdot p_0$. Because $\lambda_0 < 1$, this upper bound tends to 0 exponentially as $n$ increases. $\blacksquare$

If we are prepared to make one additional rationality assumption, we can also get the stronger result that, in light of $H^*$, no rational agent should give $H$ any credence.

iv) *Anticipation*. If the agent knows in advance that their credence in $H$ upon learning the outcome of the draw would not exceed $x$, their credence in $H$ should not exceed $x$.

The reason for separating out this assumption from the others is that, in the context of countable additivity, reflection principles are very controversial. Still, this particular claim seems plausible, and does not require a commitment to any more general reflection principle.

**Result 2**. If the agent satisfies conditions (i-iv), then their credence in $H$ prior to any draws must already be equal to 0.

> *Proof.* Let $\text{Cr}^+$ be the agent's rational credence conditioned on $H \vee H^*$ posterior to learning the outcome of a draw. The agent knows that $\bigvee_k (\text{Cr}^+(H) = \text{Cr}(H | E_k))$, which entails that $\text{Cr}^+(H) \le \lambda \cdot p$ by the lemma. So by *Anticipation*, $\text{Cr}(H) = p \le \lambda \cdot p$. So $p = 0$ (for if $p$ were positive, then $p > \lambda \cdot p$ since $\lambda < 1$). $\blacksquare$